# Dyck Numbers, II.
# Triplets and Rooted Trees in OEIS A036991

Gennady Eremin

ergenns@gmail.com

November 2, 2022

**Abstract.** Dyck paths are among the most heavily studied Catalan families. This paper is a continuation of [2]. In the paper we are dealing with the numbering of Dyck paths, the terms of the OEIS sequence A036991 or Dyck numbers. We consider triplets of terms of the form (t-4, t-2, t) (t is the senior term); triplets cover 80% of A036991. Triplets include all Mersenne numbers, with each Mersenne number being a senior term in some triplet. Triples are distributed over A036991 ranges; both the length of the ranges and the number of triplets in the range are counted by the terms of the OEIS sequence A001405. In addition to triplets, there are many lone terms in the ranges, which are counted by the terms of the OEIS sequence A116385. Each lone term (there are an infinite number of them) is the root of an infinite ternary tree of triplets. As a result, the sequence A036991 is a forest of such directed trees. We test the twin prime conjecture on infinite ternary trees of A036991. We consider the infinity theorem for pairs of prime numbers that differ by no more than 4.

*Keywords*: Dyck path, Dyck numbers, OEIS, Mersenne numbers, triplet, root tree, unary tree, ternary tree, forest of root trees, twin prime conjecture.

## 1 Introduction

Dyck paths are among the most heavily studied Catalan families [1]. This article is a continuation of [2]. In the paper we are dealing with the numbering of Dyck paths, the terms of the OEIS sequence A036991 [3] or *Dyck numbers*. A Dyck path is a directed lattice path from the origin giving up-steps U = (1, 1) and down-steps D = (1, −1) in the upper half plane and are conditioned to end in the $x$-axis. In A036991 terms, the following binary coding is adopted: an up-step is coded by zero, and a down-step is coded by 1. As a result, we get sufficiently compact decimal integers. Let's show the first A036991 terms, in red we have highlighted the well-known Mersenne numbers (see the OEIS A000225):

(1) D = 1, 3, 5, 7, 11, 13, 15, 19, 21, 23, 27, 29, 31, 39, 43, 45, 47, 51, 53, 55, 59 …

The same terms are shown below in binary expansion, in which Dyck paths are clearly traced (see the OEIS A350346). As is customary in integers, leading zeros (initial up-steps) are omitted in binary codes:

(2) B = 1, 11, 101, 111, 1011, 1101, 1111, 10011, 10101, 10111, 11011, 11101 …

In sequences (1) and (2), the terms are arranged in order of increasing length of binary code. A group of terms of the same length is called a *range*. For example, in the 4th range (4-range) there are three binary terms of length 4: $11 = 1011_2$, $13 = 1101_2$ and $15 = 1111_2$.

In each range, the last binary term is a *repunits* (repeating units), a term from the OEIS A002275:  0, 1, 11, 111, 1111, 11111 …

In both sequences, we have removed the initial null term (corresponding to an empty Dyck path) because we are working with real paths. And many sources do as well.

In sequences D and B, we will index the terms from one, that is, $d_1 = 1$, $d_2 = 3$, $b_1 = 1$, $b_2 = 11$, and so on.  Mersenne numbers divide the sequences (1) and (2) into ranges:  the *n*-th Mersenne number ends the *n*-range, i.e. the range in which all terms have length *n* in binary code. The sizes of the ranges correspond to the terms in the OEIS sequence A001405. Here is the beginning of A001405 (terms are numbered from 0):

(3)  1, 1, 2, 3, 6, 10, 20, 35, 70, 126, 252, 462, 924, 1716, 3432, 6435, 12870 …

# 2 Unary trees

A well-known recursive formula for Mersenne numbers  $M_{n+1} = 2M_n + 1$, and this can be easily verified in sequence (1). And interestingly, a similar formula is also true for any A036991 term. This is simply proved by binary expansion: if we add zero to the binary code of an arbitrary Dyck number (equivalent to multiplying by 2 in decimal form) and replace the final 0 with 1 (adding 1 in decimal form), then we get a new Dyck number from the next range. Let’s formulate the corresponding statement.

**Proposition 1.** Let *d* be an A036991 term from the *n*-range. Then

$$d' = 2d + 1 \qquad (4)$$

is the A036991 term from the (*n*+1)-range. Let's say the term *d* generates term *d’*.

Usually, the term *d’* is called the *child* of *d*, and term *d* is the *parent* of *d’*. In this case, each term has one child, but not each term has a parent.  It is not difficult to see that the inverse formula is not always true, since the removal of a last binary unit can break the *dynamics of binary code* in A036991/A350346 terms (*in each binary suffix, the number of 0’s does not exceed the number of 1’s*).  For example,  for the term 19,  the number (19 – 1)/2 = 9 is not an A036991 term.  Let’s trace the connections between terms of different ranges using formula (4); we will do it below in tabular form.

| # | Terms of OEIS A036991 | | | | | | | | | | | | | | | | | | | |
|---|---|---|---|---|---|---|---|---|---|---|---|---|---|---|---|---|---|---|---|---|
| 1 | | | | | | | | | | | | | | | | | | | | 1 |
| 2 | | | | | | | | | | | | | | | | | | | | 3 |
| 3 | | | | | | | | | | 5 | | | | | | | | | | 7 |
| 4 | | | | | | | | | | 11 | | | | 13 | | | | | | 15 |
| 5 | | | | 19 | | | 21 | | | 23 | | | | 27 | | | 29 | | | 31 |
| 6 | | | | 39 | | | 43 | 45 | | 47 | 51 | 53 | | 55 | | | 59 | 61 | | 63 |
| 7 | 71 | 75 | 77 | 79 | 83 | 85 | 87 | 91 | 93 | 95 | 103 | 107 | 109 | 111 | 115 | 117 | 119 | 123 | 125 | 127 |

Table 1. Unary trees with top roots.

In Table 1, the first column shows the range numbers. Each column of terms can be viewed as a directed unary tree with a top *root* (starting node). The right column, the ini-tial unary tree rooted at 1, is made up of the Mersenne numbers with which the ranges end. The second unary tree rooted at 5 starts in the 3-range. Note that term 5 has no parent, and this is the case for every root of this unary tree.

Term 5 is the successor of the 2nd Mersenne number: $DS(M_2) = M_2 + M_1 + 1 = 5$, where the $DS$-function is the function of determining the successor of the A036991 term. The next 4-range begins with term $DS(M_3) = M_3 + M_2 + 1 = 11$, and so on. Recall the general formula for Mersenne numbers [2]

$$DS(M_n) = M_n + M_m + 1, \; m = \lceil n/2 \rceil .$$

The number of roots grows from range to range. For example, in the 7-range, 10 new trees appear with the following roots: 71, 75, 77, 83, 85, 93, 109, 115, 117, and 125.

Obviously, the number of new roots in each range is determined by the difference between the length of the current range and the length of the previous range (see (3)). Such numbers correspond to the terms in the OEIS sequence A037952. Here is the beginning of A037952 (terms are numbered from 0):

0, 1, 1, 3, 4, 10, 15, 35, 56, 126, 210, 462, 792, 1716, 3003, 6435, 11440 …

We also need a sequence of gaps (intervals) between Dyck numbers: $g_n = DS(d_n) - d_n$. The sequence of gaps between Dyck numbers is not yet in the OEIS [3], so we will show significantly more elements (the gap between the Mersenne number and the subsequent term is indicated in red):

G = 2, 2, 2, 4, 2, 2, 4, 2, 2, 4, 2, 2, 8, 4, 2, 2, 4, 2, 2, 4, 2, 2, 8, 4, 2, 2, 4, 2, 2, 4, 2, 2, 8, 4, 2, 2, 4, 2, 2, 4, 2, 2, 16, 8, 4, 2, 2, 8, 4, 2, 2, 4, 2, 2, 4, 2, 2, 8, 4, 2, 2, 4, 2, 2, 4, 2, 2, 8, 4, 2, 2, 4, 2, 2, 4, 2, 2, 16, 8, 4, 2, 2, 8, 4, 2, 2, 4, 2, 2, 4, 2, 2, 8, 4, 2, 2, 4, 2, 2, 4, 2, 2, 8, 4, 2, 2, 4, 2, 2, 4, 2, 2, 16, 8, 4, 2, 2, 8, 4, 2, 2, 4, 2, 2, 4, 2, 2, 8, 4, 2, 2, 4, 2, 2, 4, 2, 2, 8, 4, 2, 2, 4, 2, 2, 4, 2, 2, 32, 16, 8, 4, 2, 2, 16, 8, 4, 2, 2, 8, 4, 2, 2, 4, 2, 2, 4, 2, 2, 16, 8, 4, 2 …

All terms in (2) are odd (binary codes end in 1, otherwise the dynamics in A036991 terms is violated), so all numbers in G are even. In (2), binary suffixes of length 2 (*2-suffixes*) have the form 01 or 11 = 01 +10. Accordingly, gap 2 is often found in G. Binary suffixes of length 3 (*3-suffixes*) can only be of three types: 011, 101 = 011 + 10 and 111 = 101 + 10 (otherwise the dynamics of binary codes is violated). The corresponding gaps are two adjacent twos, common in G. Note there are no single twos.

The interval between adjacent Dyck numbers is powers of two. Let's formulate the corresponding statement (the proof is given below in Section 4.2).

**Theorem 2.** In A036991, partial differences are powers of 2.

There are many pairs of twos in the gap-sequence $G$; each pair of twos corresponds to three neighboring odd numbers, *triplet*. Triplets cover 80% of A036991 (looked at a million terms). Here are the first A036991 triplets (this sequence is not yet in the OEIS):

(1, 3, 5), (3, 5, 7), (11, 13, 15), (19, 21, 23), (27, 29, 31), (43, 45, 47) …

Each Mersenne number is included in a triplet. Everything is defined by a binary 3-suffix. There are no zeros in the binary expansion of the Mersenne number. If we subtract 2 from the $n$-th Mersenne number, $n > 2$, we get the A036991 term with the binary suffix 111 – 10 = 101. If we subtract 2 again, we get a new term with the binary suffix 101 – 10 = 011. So, the following proposition is obvious.

**Proposition 3.** $n$-th Mersenne number, $n > 2$, completes a A036991 triplet of the form $(M_n - 4, M_n - 2, M_n)$.

There are an infinite number of Mersenne numbers, so the following proposition is true.

**Proposition 4.** There are an infinite number of triplets in the OEIS A036991.

Let's remember the twin primes (the OEIS sequence A001097); there are a lot of twin primes in A036991 triplets. If we find at least one triplet with the twin primes in each range (we will even be satisfied if the primes are placed on the triplet edges), then we can prove the hypothesis of the infinity of twin primes.

In paper [2], we worked with the successor function of the Dyck number, DS-function. Let's try to get something similar for A036991 triplets.

# 3 Triplet generation

**3.1.** Let us count the number of triplets in the ranges of (1). Recall that in A036991 ranges, the terms are grouped by binary code length. For example, the 1st and 2nd ranges contain one term each, and in the 3-range there are two terms 5 and 7. Accordingly, there are no triplets. But in the next 4-range there are three terms 11, 13 and 15, which form a triplet, and so on. Below in Table 2 we have shown the data of the first 8 ranges.

| Range | Number of terms | Terms and triplets | Number of triplets |
|---|---|---|---|
| 1 | 1 | 1 | 0 |
| 2 | 1 | 3 | 0 |
| 3 | 2 | 5, 7 | 0 |
| 4 | 3 | (11, 13, 15) | 1 |
| 5 | 6 | (19, 21, 23), (27, 29, 31) | 2 |
| 6 | 10 | 39, (43, 45, 47), (51, 53, 55), (59, 61, 63) | 3 |
| 7 | 20 | 71, (75, 77, 79), (83, 85, 87), (91, 93, 95),<br>103, (107, 109, 111), (115, 117, 119), (123, 125, 127) | 6 |
| 8 | 35 | 143, 151, (155,157,159),<br>167, (171, 173, 175), (179, 181, 183), (187, 189, 191),<br>199, (203, 205, 207), (211, 213, 215), (219, 221, 223),<br>231, (235, 237, 239), (243, 245, 247), (251, 253, 255) | 10 |

Table 2. Distribution by ranges of single terms and triplets.

Some patterns are obvious. The range sizes for terms and triplets are the same with a shift of two ranges. For example, we have 10 terms in 6-range, and the same number of triplets is in 8-range. Recall that the range sizes correspond to the A001405 terms, see (3). To match Table 2, we will index the numbers in (3) from 1. The first 16 ranges are easy to check using the A036991 *b*-file. Let's look at a small example.

**Example 5**. In the b-file of A036991 (see Table in Section LINKS), let's extract the extreme terms from the 16-range (see *DS*-formula in [2]):

the first term $DS(M_{15}) = M_{15} + M_8 + 1 = 32767 + 255 + 1 = 33023 = \text{A036991}(7062)$,
the last term $M_{16} = 65535 = \text{A036991}(13496)$.

As a result we get the length of the 16-range 13496 – 7062 + 1 = 6435 = A001405(16). It remains for us to show that in this range the number of triplets is A001405(14) = 1716.

**3.2.** We are looking for connections between terms that are separated by two ranges. Let's go back to formula (4) $d' = 2d + 1$, and take one more step, moving to the next range:

(5) $$d'' = 2d' + 1 = 2\times(2d + 1) + 1 = 4d + 3.$$

Let's analyze the resulting formula in binary expansion. The original term $d$, as usual, ends with the binary 1-suffix 1. Multiplying by $4 = 100_2$, we just add two zeros at the end of the binary code, so the new number ends with the binary 3-suffix 100. And finally, by adding $3 = 11_2$, we get the final binary 3-suffix 111 in $d''$. Thus, we get the last term of some triplet, as a result we get access to two additional terms $d'' - 2$ and $d'' - 4$. Let's say that the given term $d$ from $n$-range using formula (5) *generated* in the $(n+2)$-range a triplet of the form

(6) $$(d''- 4, d''- 2, d'') = (4d - 1, 4d + 1, 4d + 3).$$

**Proposition 6.** The generation operation is unique; the resulting triplets do not have common terms (they do not intersect).

*Proof*. Let's perform the generation procedure for all terms of the triplet $(d, d + 2, d + 4)$. In addition to (6), we need to process $d + 2$ and $d + 4$. As a result we get two more triplets:

$d + 2 \Rightarrow (4d + 7, 4d + 9, 4d + 11)$,
$d + 4 \Rightarrow (4d + 15, 4d + 17, 4d + 19)$.

As we see, triplet (6) and both new triplets have no common terms. □

**Lemma 7.** Let $D_n$ be the $n$-range in A036991 and $T_n$ be the triplet set in the $n$-range. Then

$$|T_n| = |D_{n-2}| .$$

*Proof*. The A036991 terms do not repeat and hence according to (6) the triplets also do not repeat. Conversely, it is easy to show that each triplet corresponds to a certain term. It is enough to remove the last three digits in the binary expansion of each triplet term, and then add a binary unit; as a result, we will get the same term. In decimal form, we get the corresponding formula that is inverse to (5):

(7) $$d = 2\times(d''//8) + 1.$$

We borrowed the integer division operation '//' from the Python language (the fractional part is discarded). It is not difficult to see that formula (7) is true if we replace $d''$ with $d''- 2$ or $d''- 4$ (recall that the term $d''$ has the binary 3-suffix 111).

As a result, we obtain a one-to-one correspondence (bijection) between elements of set $T_n$ and set $D_{n-2}$. □

Below is another small example.

**Example 8**. Let's run the generation procedure (6) step by step, starting with the terms 3, 5, 7 (this first A036991 triplet we got earlier from 1).

3 => (11, 13, 15) => (43, 45, 47), (51, 53, 55), (59, 61, 63) => (171, 173, 175) …
5 => (19, 21, 23) => (75, 77, 79), (83, 85, 87), (91, 93, 95) => (299, 301, 303) …
7 => (27, 29, 31) => (107,109,111), (115,117,119), (123,125,127) => (427,429,431) …

With the help of such generation, in three chains we got all the triplets of the initial ranges. In the chains, at each step, the number of triplets is trebled, because in each triplet, each term generates a new triplet (three children). Obviously, chains can be continued indefinitely. In the first chain we work with even ranges: 3 represents the 2-range and generates one new triplet in the 4-range, term 11 from the 4-range generates the first triplet in the 6-range (total 3 new triplets), 43 starts the first triplet from the 6-range (total 9 new triplets), 171 starts the first triplet from the 8-range (total 27 new triplets), and so on. In the other two chains, triplets of odd ranges are generated.

# 4 Ternary root trees

**4.1.** The term 1 generates the first triplet (3, 5, 7), and it is the only triplet with terms from different ranges: term 3 is from 2-range, and terms 5, 7 are from 3-range. Therefore term 3 generates triplets in even ranges, while terms 5 and 7 spawn triplets of odd ranges. As a result, we get a vast infinite ternary tree with the root 1. In other words, we get a directed root tree all edges of which point away from the root 1. Each node of such a tree has three children, i.e. we get a ternary directed tree.

**Definition 9.** Let's call the directed tree with the root 1 a *base ternary tree*.

There are other infinite root trees, and there are infinitely many of them, but the base tree is the most extensive (within each finite range). For example, the triplets of the base tree are placed on both even and odd ranges. This is not the case for other root trees. Another important fact is that all Mersenne numbers are contained in triplets of the base tree.

**Definition 10.** A term that is not included in a triplet is called a *lone term*.

The interval between a lone term and neighboring terms is greater than 2; it's 4, 8, 16, and so on. Note that for a lone term, the left interval always exceeds the right one. In A036991, there are many lone terms, and each lone term generates unique triplet. Let's look at the corresponding example.

**Example 11**. In A036991 in the 6-range, there is a lone term 39 (see Table 2), which generate triplets in even ranges 8, 10, 12, and so on. In the next 7-range, we already have two lone terms 71 and 103, which generate triplets in odd ranges 9, 11, 13, and so on.

39 => (155, 157, 159) => (619, 621, 623), (627, 629, 631), (635, 637, 639)

=> (2475, 2477, 2479), (2483, 2485, 2487), (2491, 2493, 2495), (2507 …

71 => (283, 285, 287) => (1131, 1133, 1135), (1139, 1141, 1143), (1147, 1149, 1151)
=> (4523, 4525, 4527), (4531, 4533, 4535), (4539, 4541, 4543), (4555 …

103 => (411, 413, 415) => (1643, 1645, 1647), (1651, 1653, 1655), (1659, 1661, 1663)
= > (6571, 6573, 6575), (6579, 6581, 6583), (6587, 6589, 6591), (6603 …

Thus, in the 6-range, in addition to the base tree, the lone term 39 starts another infinite ternary tree. In the next 7-range, two infinite ternary trees with root 71 and root 103 are launched at once.

**Proposition 12.** In A036991, each lone term from an even (odd) range generates in even (odd) ranges a unique infinite directed ternary tree, in which such a term is a root.

We can say that the transition from (4) to (5) allowed us to move from unary trees to ternary ones.

**4.2.** A natural question arises as to how many lone terms and how they appear in A036991. Table 3 below shows the distribution of lone terms by ranges.

| Range | Number of lone terms | Lone terms (roots of ternary trees) |
|---|---|---|
| 6 | 1 | 39 |
| 7 | 2 | 71, 103 |
| 8 | 5 | 143, 151, 167, 199, 231 |
| 9 | 10 | 271, 279, 295, 327, 359, 399, 407, 423, 455, 487 |
| 10 | 21 | 543, 559, 567, 591, 599, 615, 655, 663, 679, 711, 743, 783, 791, 807, 839, 871, 911, 919, 935, 967, 999 |
| 11 | 42 | 1055, 1071, 1079, 1103, 1111, 1127, 1167, 1175, 1191, 1223, 1255, 1295, 1303, 1319, 1351, 1383, 1423, 1431, 1447, 1479, 1511, 1567, 1583, 1591, 1615, 1623, 1639, 1679, 1687, 1703, 1735, 1767, 1807, 1815, 1831, 1863, 1895, 1935, 1943, 1959, 1991, 2023 |

Table 3. Lone terms in A036991 (roots of ternary trees).

Table 3 shows 1+2+5+10+21+42 = 81 roots of ternary trees. The sequence of roots is infinite; it is not yet in the OEIS. It is easy to count the number of lone terms in any range of A036991 if we return to A001405 (3). To get the desired value, it is enough to subtract three times the number of triplets from the length of the range. For example, in the 8th range we get $35 - 3\times10 = 5$ lone terms (see Table 2). The resulting expression gives us the terms of the OEIS sequence A116385. Here is the beginning of A116385:

0, 0, 1, 2, 5, 10, 21, 42, 84, 168, 330, 660, 1287, 2574, 5005, 10010, 19448 …

In the OEIS, the A001405 terms and the A116385 terms are indexed from zero, so the resulting dependence takes the form:

$$\text{(8)} \qquad \text{A116385}(n) = \text{A001405}(n+3) - 3\times\text{A001405}(n+1), \quad n \geq 0.$$

In the OEIS sequence A116385, there is no formula (8) yet.

The Appendix A shows a list of the roots of more than three hundred ternary trees.

Let's describe how lone terms are formed. In A036991, a term with a large binary repunit suffix (repeated units) [4] is often followed by a lone term with a repunit suffix twice as short. For example, this is true for Mersenne numbers, in the binary expansion of which there are no zeros. For the $n$-th Mersenne number, and many other terms $d$ with a repunit $n$-suffix, the next term is determined by the following formula (see [2]):

$$\text{(9)} \qquad d' = DS(d) = d + 2^m, \text{ where } m = \lceil n/2 \rceil.$$

Recall the $DS$-function is the function of determining the successor of the A036991 term.

Additionally, we note that lone terms are formed only in the case $n > 4$. Accordingly, the subsequent terms can also be lone. Here are these lonely terms:

$$\text{(10)} \qquad d'' = DS(d') = d' + 2^{m-1}, \quad d''' = DS(d'') = d'' + 2^{m-2}, \text{ and so on.}$$

*Proof of Theorem 2*. Equality (9) and sequence (10) confirm Theorem 2. Another way to obtain a term successor is to exchange the high digit of the repunit suffix and the left adjacent zero in binary expansion. Such an exchange does not violate the dynamics of the binary code in A036991 terms. For a repunit $n$-suffix, we get a successor increment of the form $2^n - 2^{n-1} = 2^{n-1}$ (recall that in integers the digits are numbered from right to left, starting from zero). □

The Appendix B shows the triplets of the initial levels (ranges) of a ternary tree with root 39.

**4.3.** Regarding sequence (10), in the case of Mersenne numbers, $d = M_n$, as $n$ increases, we get an increasing (in the limit infinite) sequence of lone terms in the $(n+1)$-th range. Each lone term is the root of an infinite tree of triplets. You can also add that we have an infinite number of ranges (each range ends with a Mersenne number, the number of which is infinite). We can say that we have infinity in three dimensions, or *infinity cubed*. Let us formulate the corresponding theorem.

**Theorem 13.** The OEIS A036991 is a *forest* of infinite directed ternary trees. The number of trees in such a forest is infinite.

# 5 Dyck primes, twin prime conjecture

**5.1.** There are a lot of prime numbers in the OEIS A036991, up to 30% (a million terms have been checked). The situation is similar with primes (see the OEIS A000040): in the first million A000040 terms, 304208 numbers are A036991 terms.

**Definition 14.** Let's called the prime terms of A036991 *Dyck primes*.

Here are the first Dyck primes (see the OEIS A350577):

3, 5, 7, 11, 13, 19, 23, 29, 31, 43, 47, 53, 59, 61, 71, 79, 83, 103, 107, 109 …

As we can see, among the Dyck primes there are many twin primes (see the OEIS A001097): 3, 5, 7, 11, 13, 29, 31, 59, 61, 107, 109, 179, 181, 347, 349, 461, 463 …

As you know, there are infinitely many prime numbers (proved by Euclid 300 BC). It would also be good for us to prove the infinity of the sequence of Dyck primes (in A350577, the author has so far formulated the corresponding conjecture).

As for the twin primes, everything is much more complicated. It is currently unknown whether there are infinitely many twin primes; the conjecture of twin primes has not yet been proven. For the past 10 years, mathematicians have been actively working on this problem. So far, it has been proved that there are infinitely many pairs of primes that differ by no more than 246 [5]. In this paper, it is proposed to work with the conjecture of twin primes using infinite root trees from the OEIS sequence A036991. Recall that the number of such ternary trees is also infinite.

A036991 terms are distributed in sequence by ranges (or strings), according to the length of binary codes. Each $n$-th string ends with the $n$-th Mersenne number, whose binary expansion has no zeros (see Table 2). And the same string begins with the successor of the ($n$-1)-th Mersenne number. The number of Mersenne numbers is infinite, so if we prove that there is at least one prime in each string (the length of the string increases with the growth of $n$), then we get an infinity of Dyck primes.

**Definition 15.** A triplet that includes two prime numbers is called a *prime triplet*.

Among three consecutive odd numbers, one is always a multiple of three. Therefore, in a prime triplet, in addition to two primes, the third term is a multiple of three. Obviously, a prime triplet contains either a pair of twin primes or two prime numbers that differ by 4. In both cases, if we prove the infinity of prime triplets, we will greatly improve the last result 246, since 4 is much smaller.

**5.2.** The Appendix C shows the prime triplets of initial ranges. For compactness, terms in triplets are connected by a "/", and composite terms, multiples of 3, are zeroed (shown in red). As you can see, there are no triplets in the first three short ranges, but further, as the length of the ranges increases, the number of prime triplets constantly increases (and certainly does not decrease, and this is very important!). Here is the number of prime triplets in the first 22 ranges:

0, 0, 0, 1, 2, 2, 1, 1, 5, 5, 10, 17, 17, 48, 67, 111, 207, 349, 599, 1102, 1879, 3290.

The last tested triplet is 0/4192757/4192759; in the binary expansion, the length of these terms is 22.

**5.3.** The last Appendix D shows prime triplets in a ternary tree rooted at 39. Again, terms in triplets are separated by a "/", and terms that are multiples of 3 are set to zero. The first prime triplets appear in the 12-range (in such terms the length of the binary code is 12). And in this case, with increasing range length, the number of prime triplets is constantly growing.

Let me remind you that to test the twin prime conjecture, it is enough to show that in such a ternary tree there is at least one prime triplet. There are an infinite number of such ternary trees in the OEIS sequence A036991. I think we can formulate a corresponding theorem.

**Theorem 16.** There are infinitely many pairs of prime numbers that differ by no more than four.

# 6 Conclusions and future work

We can say that the OEIS sequence A036991 is an empire of THREE. Firstly, the vast majority of terms are combined into *triplets* (three adjacent odd numbers). Secondly, triplets are grouped into additional triples, and as a result we get *Nines* (*three squared*); the interval between adjacent triplets is 4. Thirdly, Nines are also grouped into triples, and in such a group there are already *twenty-seven* terms (*three cubed*). Adjacent Nines are separated by lone terms, which are the roots of directed ternary trees. For example, the first lone term 39 appears between Nine 11…31 (we show extreme terms) and Nine 43…63. Next comes the second lone term 71, followed by the third Nine 75…95, and so on.

This triple enlargement procedure can be continued further. And there seems to be no end to it. And this is a topic for future work. The author will be grateful to readers who are interested in this subject and will take an active part in this project.

**Acknowledgements.** The author would like to thank Jörg Arndt (Technische Hochschule Nürnberg), Alois P. Heinz (University of Applied Sciences Heilbronn), and Olena G. Kachko (Kharkiv National University of Radio Electronics) for discussion of the considered integer sequences.

# References

[1] R. Stanley. *Catalan numbers*. Cambridge University Press, Cambridge, 2015.

[2] Gennady Eremin. *Dyck numbers, I. Successor Function*, 2022. arXiv:2210.00744

[3] Neil J. A. Sloane, Ed. *The On-Line Encyclopedia of Integer Sequences*. https://oeis.org/, 2022.

[4] Wolfram MathWorld, Repunit, 2022. https://mathworld.wolfram.com/Repunit.html

[5] Bounded gaps between primes. Polymath. Retrieved 2014-03-27
https://asone.ai/polymath/index.php?title=Bounded_gaps_between_primes

Mentions the OEIS sequences A000040, A000225, A001097, A001405, A002275, A005408, A036991, A037952, A116385, A350346, A350577.

Gzhel State University, Moscow, 140155, Russia
http://www.art-gzhel.ru/

## A Lone terms (roots of ternary trees)

**Range 6**: 39.

**Range 7**: 71, 103.

**Range 8**: 143, 151, 167, 199, 231.

**Range 9**: 271, 279, 295, 327, 359, 399, 407, 423, 455, 487.

**Range 10**: 543, 559, 567, 591, 599, 615, 655, 663, 679, 711, 743, 783, 791, 807, 839, 871, 911, 919, 935, 967, 999.

**Range 11**: 1055, 1071, 1079, 1103, 1111, 1127, 1167, 1175, 1191, 1223, 1255, 1295, 1303, 1319, 1351, 1383, 1423, 1431, 1447, 1479, 1511, 1567, 1583, 1591, 1615, 1623, 1639, 1679, 1687, 1703, 1735, 1767, 1807, 1815, 1831, 1863, 1895, 1935, 1943, 1959, 1991, 2023.

**Range 12**: 2111, 2143, 2159, 2167, 2207, 2223, 2231, 2255, 2263, 2279, 2335, 2351, 2359, 2383, 2391, 2407, 2447, 2455, 2471, 2503, 2535, 2591, 2607, 2615, 2639, 2647, 2663, 2703, 2711, 2727, 2759, 2791, 2831, 2839, 2855, 2887, 2919, 2959, 2967, 2983, 3015, 3047, 3103, 3119, 3127, 3151, 3159, 3175, 3215, 3223, 3239, 3271, 3303, 3343, 3351, 3367, 3399, 3431, 3471, 3479, 3495, 3527, 3559, 3615, 3631, 3639, 3663, 3671, 3687, 3727, 3735, 3751, 3783, 3815, 3855, 3863, 3879, 3911, 3943, 3983, 3991, 4007, 4039, 4071.

**Range 13**: 4159, 4191, 4207, 4215, 4255, 4271, 4279, 4303, 4311, 4327, 4383, 4399, 4407, 4431, 4439, 4455, 4495, 4503, 4519, 4551, 4583, 4639, 4655, 4663, 4687, 4695, 4711, 4751, 4759, 4775, 4807, 4839, 4879, 4887, 4903, 4935, 4967, 5007, 5015, 5031, 5063, 5095, 5151, 5167, 5175, 5199, 5207, 5223, 5263, 5271, 5287, 5319, 5351, 5391, 5399, 5415, 5447, 5479, 5519, 5527, 5543, 5575, 5607, 5663, 5679, 5687, 5711, 5719, 5735, 5775, 5783, 5799, 5831, 5863, 5903, 5911, 5927, 5959, 5991, 6031, 6039, 6055, 6087, 6119, 6207, 6239, 6255, 6263, 6303, 6319, 6327, 6351, 6359, 6375, 6431, 6447, 6455, 6479, 6487, 6503, 6543, 6551, 6567, 6599, 6631, 6687, 6703, 6711, 6735, 6743, 6759, 6799, 6807, 6823, 6855, 6887, 6927, 6935, 6951, 6983, 7015, 7055, 7063, 7079, 7111, 7143, 7199, 7215, 7223, 7247, 7255, 7271, 7311, 7319, 7335, 7367, 7399, 7439, 7447, 7463, 7495, 7527, 7567, 7575, 7591, 7623, 7655, 7711, 7727, 7735, 7759, 7767, 7783, 7823, 7831, 7847, 7879, 7911, 7951, 7959, 7975, 8007, 8039, 8079, 8087, 8103, 8135, 8167 (total 334 roots).

## B Ternary tree with root 39 (even ranges of A036991).

**Range 6** (1 node): root 39.

**Range 8** (3 nodes): 155, 157, 159.

**Range 10** (9 nodes): 619, 621, 623, 627, 629, 631, 635, 637, 639.

**Range 12** (27 nodes): 2475, 2477, 2479, 2483, 2485, 2487, 2491, 2493, 2495, 2507, 2509, 2511, 2515, 2517, 2519, 2523, 2525, 2527, 2539, 2541, 2543, 2547, 2549, 2551, 2555, 2557, 2559.

**Range 14** (81 nodes): 9899, 9901, 9903, 9907, 9909, 9911, 9915, 9917, 9919, 9931, 9933, 9935, 9939, 9941, 9943, 9947, 9949, 9951, 9963, 9965, 9967, 9971, 9973, 9975, 9979, 9981, 9983, 10027, 10029, 10031, 10035, 10037, 10039, 10043, 10045, 10047, 10059, 10061, 10063, 10067, 10069, 10071, 10075, 10077, 10079, 10091, 10093, 10095, 10099, 10101, 10103, 10107, 10109, 10111, 10155, 10157, 10159, 10163, 10165, 10167, 10171, 10173, 10175, 10187, 10189, 10191, 10195, 10197, 10199, 10203, 10205, 10207, 10219, 10221, 10223, 10227, 10229, 10231, 10235, 10237, 10239.

## C Prime triplets in the initial ranges,

composite terms, multiples of 3, are set to red zero

**Range 4** (3 terms): 11/13/0 (1 prime triplet).

**Range 5** (6 terms): 19/0/23, 0/29/31 (2 prime triplets).

**Range 6** (10 terms): 43/0/47, 59/61/0 (2 prime triplets).

**Range 7** (20 terms): 107/109/0 (1 prime triplet).

**Range 8** (35 terms): 179/181/0 (1 prime triplet).

**Range 9** (70 terms): 307/0/311, 347/349/0, 379/0/383, 0/461/463, 499/0/503 (5 prime triplets).

**Range 10** (126 terms): 0/821/823, 827/829/0, 859/0/863, 883/0/887, 1019/1021/0 (5 prime triplets).

**Range 11** (252 terms): 0/1229/1231, 0/1277/1279, 1451/1453/0, 1483/0/1487, 1787/1789/0, 1867/0/1871, 0/1877/1879, 0/1949/1951, 0/1997/1999, 2027/2029/0 (10 prime triplets).

**Range 12** (462 terms): 0/2237/2239, 2267/2269/0, 2539/0/2543, 0/2549/2551, 2683/0/2687, 3019/0/3023, 3163/0/3167, 3187/0/3191, 3251/3253/0, 3371/3373/0, 0/3389/3391, 3539/3541/0, 0/3581/3583, 0/3821/3823, 0/3917/3919, 4019/4021/0, 4091/4093/0 (17 prime triplets).

**Range 13** (924 terms): 4787/4789/0, 0/5021/5023, 5099/5101/0, 5227/0/5231, 0/5501/5503, 0/5741/5743, 5867/5869/0, 6043/0/6047, 6131/6133/0, 0/6269/6271,

6779/6781/0, 6827/6829/0, 0/6869/6871, 6907/0/6911, 0/7349/7351, 7547/7549/0, 7603/0/7607 (17 prime triplets).

**Range 14** (1716 terms): 8443/0/8447, 0/9341/9343, 0/9437/9439, 0/9461/9463, 0/9629/9631, 0/9677/9679, 9787/0/9791, 9883/0/9887, 0/10037/10039, 10067/10069/0, 10091/10093/0, 10099/0/10103, 10139/10141/0, 10427/10429/0, 10459/0/10463, 0/10709/10711, 10859/10861/0, 0/11069/11071, 11083/0/11087, 0/11117/11119, 11443/0/11447, 11467/0/11471, 0/11549/11551, 11699/11701/0, 12107/12109/0, 12251/12253/0, 12539/12541/0, 12907/0/12911, 0/12917/12919, 12979/0/12983, 13003/0/13007, 13099/0/13103, 13147/0/13151, 13691/13693/0, 0/13757/13759, 13931/13933/0, 0/13997/13999, 14323/0/14327, 14779/0/14783, 14827/0/14831, 15259/0/15263, 0/15581/15583, 15643/0/15647, 15667/0/15671, 15731/15733/0, 15787/0/15791, 0/16061/16063, 16187/16189/0 (48 prime triplets).

**Range 15** (3432 terms): 0/16829/16831, 17387/17389/0, 0/17597/17599, 0/17789/17791, 0/17837/17839, 0/17909/17911, 18043/0/18047, 18131/18133/0, 18251/18253/0, 0/19181/19183, 19379/19381/0, 19387/0/19391, 0/19421/19423, 20147/20149/0, 0/20477/20479, 0/20717/20719, 0/20981/20983, 21491/21493/0, 21499/0/21503, 21611/21613/0, 0/22109/22111, 23027/23029/0, 23291/23293/0, 0/23669/23671, 0/23741/23743, 24019/0/24023, 24179/24181/0, 24371/24373/0, 25307/25309/0, 25339/0/25343, 0/25469/25471, 25579/0/25583, 26107/0/26111, 0/26861/26863, 27059/27061/0, 0/27581/27583, 27739/0/27743, 27763/0/27767, 28027/0/28031, 0/28109/28111, 0/28277/28279, 0/28349/28351, 28571/28573/0, 28603/0/28607, 28619/28621/0, 28627/0/28631, 0/28661/28663, 0/29021/29023, 29387/29389/0, 30011/30013/0, 0/30269/30271, 0/30389/30391, 30427/0/30431, 30491/30493/0, 0/30557/30559, 31147/0/31151, 0/31181/31183, 31219/0/31223, 31387/0/31391, 0/31541/31543, 31723/0/31727, 32027/32029/0, 32059/0/32063, 0/32117/32119, 0/32189/32191, 32411/32413/0, 0/32717/32719 (67 prime triplets).

**Range 16** (6435 terms): 0/34301/34303, 35323/0/35327, 0/35837/35839, 36187/0/36191, 0/36341/36343, 36467/36469/0, 36523/0/36527, 36779/36781/0, 36787/0/36791, 37307/37309/0, 37691/37693/0, 37811/37813/0, 0/38237/38239, 38299/0/38303, 38459/38461/0, 38651/38653/0, 38707/0/38711, 38747/38749/0, 39227/39229/0, 0/39341/39343, 39371/39373/0, 39667/0/39671, 39883/0/39887, 40123/0/40127, 40427/40429/0, 0/40637/40639, 40819/0/40823, 0/41981/41983, 0/42221/42223, 0/42461/42463, 0/42701/42703, 42859/0/42863, 43987/0/43991, 44027/44029/0, 44203/0/44207, 44267/44269/0, 0/44381/44383, 44531/44533/0, 44699/44701/0, 45179/45181/0, 0/45821/45823, 46507/0/46511, 0/46589/46591, 46747/0/46751, 0/46829/46831, 0/47741/47743, 47947/0/47951, 48539/48541/0, 48619/0/48623, 48731/48733/0, 0/48821/48823, 0/48989/48991, 49787/49789/0, 50539/0/50543, 0/50549/50551, 50587/0/50591, 50891/50893/0, 51059/51061/0, 51131/51133/0, 0/51197/51199, 51419/51421/0, 0/51437/51439, 51827/51829/0, 0/51869/51871, 0/52181/52183, 0/52541/52543, 52859/52861/0, 53147/53149/0, 53171/53173/0, 0/53717/53719, 54011/54013/0, 0/54581/54583, 0/55661/55663, 0/55901/55903, 55931/55933/0, 0/56237/56239, 56267/56269/0, 0/56477/56479, 0/56501/56503, 56531/56533/0, 56659/0/56663, 56779/0/56783, 56891/56893/0, 57139/0/57143, 0/57269/57271, 57787/0/57791, 0/58109/58111, 0/58229/58231,

58363/0/58367, 59051/59053/0, 59219/59221/0, 0/59357/59359, 59627/59629/0, 60659/60661/0, 0/60917/60919, 61339/0/61343, 61627/0/61631, 0/62141/62143, 0/62189/62191, 62299/0/62303, 62323/0/62327, 62683/0/62687, 0/63197/63199, 0/63389/63391, 63419/63421/0, 64187/64189/0, 0/64301/64303, 0/64877/64879, 65179/0/65183, 65267/65269/0, 65323/0/65327 (111 prime triplets).

## D Prime triplets in ternary tree, root 39

(composite terms, multiples of 3, are set to red zero)

**Range 12** (27 terms): 2539/0/2543, 0/2549/2551 (2 prime triplets).

**Range 14** (81 terms): 0/10037/10039, 10067/10069/0, 10091/10093/0, 10099/0/10103 (4 prime triplets).

**Range 16** (243 terms): 39667/0/39671, 39883/0/39887, 40123/0/40127, 40427/40429/0, 0/40637/40639, 40819/0/40823 (6 prime triplets).

**Range 18** (729 terms): 158923/0/158927, 159403/0/159407, 159539/159541/0, 159667/0/159671, 160619/160621/0, 0/160637/160639, 161459/161461/0, 0/161741/161743, 161771/161773/0, 161779/0/161783, 162523/0/162527, 0/162749/162751, 163019/163021/0, 0/163061/163063, 163147/0/163151, 163307/163309/0 (16 prime triplets).

**Range 20** (2187 terms): 0/634157/634159, 634859/634861/0, 635563/0/635567, 635707/0/635711, 635891/635893/0, 636107/636109/0, 636283/0/636287, 636403/0/636407, 0/636917/636919, 0/637781/637783, 0/638717/638719, 641747/641749/0, 0/641789/641791, 642011/642013/0, 0/642797/642799, 0/642869/642871, 644051/644053/0, 644443/0/644447, 0/644597/644599, 0/644909/644911, 646099/0/646103, 646571/646573/0, 0/650477/650479, 651067/0/651071, 651179/651181/0, 0/651221/651223, 0/652541/652543, 652723/0/652727, 654163/0/654167, 654187/0/654191, 654539/654541/0, 654779/654781/0, 655219/0/655223 (33 prime triplets).

**Range 22** (6561 terms): 2534267/2534269/0, 0/2535101/2535103, 2536307/2536309/0, 2536379/2536381/0, 2536811/2536813/0, 2536907/2536909/0, 2537459/2537461/0, 2538299/2538301/0, 2538707/2538709/0, 0/2538749/2538751, 0/2539349/2539351, 2543323/0/2543327, 2544299/2544301/0, 0/2544629/2544631, 2545451/2545453/0, 2546539/0/2546543, 0/2546669/2546671, 2546899/0/2546903, 0/2546909/2546911, 0/2547029/2547031, 0/2547581/2547583, 2550763/0/2550767, 2550971/2550973/0, 2551099/0/2551103, 2551123/0/2551127, 2551499/2551501/0, 2551603/0/2551607, 0/2552621/2552623, 2552651/2552653/0, 0/2553149/2553151, 0/2554829/2554831, 0/2555261/2555263, 2567347/0/2567351, 2567531/2567533/0, 0/2568029/2568031, 2568187/0/2568191, 0/2569421/2569423, 2569939/0/2569943, 2570203/0/2570207, 2570219/2570221/0, 2571067/0/2571071, 2571731/2571733/0, 2572091/2572093/0, 2572123/0/2572127, 2575019/2575021/0, 0/2575061/2575063, 2575091/2575093/0, 2576219/2576221/0, 0/2577077/2577079, 0/2577917/2577919, 0/2578109/2578111, 0/2578349/2578351, 2579387/2579389/0, 2579803/0/2579807, 2580467/2580469/0, 0/2583389/2583391, 2583739/0/2583743, 2583859/0/2583863, 0/2585837/2585839,

0/2585981/2585983, 0/2588669/2588671, 2599627/0/2599631, 2599739/2599741/0, 2599867/0/2599871, 2600947/0/2600951, 2602331/2602333/0, 0/2602349/2602351, 0/2602877/2602879, 2604011/2604013/0, 2604731/2604733/0, 2605019/2605021/0, 0/2607989/2607991, 0/2608349/2608351, 0/2609069/2609071, 2609899/0/2609903, 2610131/2610133/0, 2610379/0/2610383, 2610611/2610613/0, 0/2610677/2610679, 0/2611157/2611159, 2612411/2612413/0, 0/2612429/2612431, 2613227/2613229/0, 0/2615981/2615983, 2616139/0/2616143, 2616667/0/2616671, 0/2616701/2616703, 2616787/0/2616791, 2617259/2617261/0, 2617267/0/2617271, 2618107/0/2618111, 0/2619389/2619391, 2620139/2620141/0, 0/2620589/2620591, 0/2620661/2620663 (95 prime triplets).